\documentclass[11pt, reqno]{amsart}
\usepackage{euscript,amscd,amsgen,amsfonts,amssymb,latexsym}

\newcommand{\eqdef}{\stackrel{\scriptscriptstyle\rm def}{=}}

\newtheorem{theorem}{Theorem}

\newtheorem{proposition}{Proposition}
\newtheorem{corollary}{Corollary}

\newtheorem{lemma}{Lemma}
\newtheorem{example}{Example}

\def\vol{{\rm vol} }
\newcommand{\beha}{\begin{enumerate}}
\newcommand{\behe}{\end{enumerate}}
\renewcommand{\epsilon}{\varepsilon}

\newcommand{\cM}{\EuScript{M}}

\newcommand{\bR}{{\mathbb R}}

\newcommand{\bZ}{{\mathbb Z}}
\newcommand{\bN}{{\mathbb N}}

\newcommand{\cL}{{\mathcal L}}
\newcommand{\cR}{{\mathcal R}}

\def\1{1\!\!1}

\def\and{\text{ and }}

         \def\P{\text{{\rm P}}}

                        \def\^{\tilde}

\def\Per{{\rm Per}}
\def\SPer{{\rm SPer}}
\def\SFix{{\rm SFix}}
\def\Fix{{\rm Fix}}

\def\card{\#}

\def\ni{\noindent}

\def\1{1\!\!1}

\DeclareMathSymbol{\varnothing}{\mathord}{AMSb}{"3F}
\renewcommand{\emptyset}{\varnothing}
\title{Topological pressure via saddle points}
\author{Katrin Gelfert}\address{Max-Planck-Institut f\"ur Physik komplexer
  Systeme, N\"othnitzer Str. 38, D-01187 Dresden \& Institut f\"ur Physik, TU
  Chemnitz, D-09107 Chemnitz, Germany}\email{gelfert@pks.mpg.de}
\urladdr{http://www.pks.mpg.de/~gelfert/}

\author{Christian Wolf}\address{Department of Mathematics, Wichita State University, Wichita, KS 67260, USA}\email{cwolf@math.wichita.edu}
\urladdr{http://www.math.wichita.edu/~cwolf/}

\begin{document}
\thanks{The research of the first author was supported by the Deutsche Forschungsgemeinschaft. 
The research of the second author was supported in part by the National Science
Foundation under Grant No. EPS-0236913 and matching support from
the State of Kansas through Kansas Technology Enterprise
Corporation.}

\begin{abstract}
Let $\Lambda$ be a compact locally maximal invariant set of a
$C^2$-diffeomorphism $f:M\to M$ on a smooth Riemannian manifold $M$.
In this paper we study the topological pressure $P_{ \rm
top}(\varphi)$ (with respect to the dynamical system $f|\Lambda$)
for a wide class of H\"older continuous potentials and analyze its
relation to dynamical, as well
as geometrical, properties of the system. We show that under a
mild nonuniform hyperbolicity assumption the  topological
pressure of $\varphi$ is entirely determined by the values of $\varphi$ on the
saddle points of $f$ in $\Lambda$.
Moreover, it is enough to consider saddle points with ``large'' Lyapunov
exponents.
We also introduce a version of the  pressure for certain
non-continuous potentials and establish several variational inequalities for
it. Finally, we deduce relations between  expansion and escape
rates and the dimension of $\Lambda$. Our results generalize several well-known results
to  certain non-uniformly hyperbolic systems.
\end{abstract}
\subjclass[2000]{}
\maketitle

\section{Introduction}

\subsection{Motivation}
In the geometric and  ergodic-theoretical aspects of the theory of
dynamical systems, the so-called thermodynamic formalism, which
was originally developed by theoretical physicists, has become a
powerful tool during the last three decades. The main object in
this theory is the topological pressure, i.e. a particular
functional on the space of observables, that encodes several
important quantities of the underlying dynamical system. For
example, pressure can be applied to obtain information about Lyapunov exponents,
dimension, multifractal spectra, natural invariant measures, etc. (see for instance~\cite{P} and the references therein).
In particular, in the case of hyperbolic systems Bowen and Ruelle
established in their pioneer works deep connections between
topological pressure and periodic points, Hausdorff dimension and
the characterization of attractors. The main purpose of this paper
is to generalize some of these results to the case of certain
non-uniformly hyperbolic systems. One key idea in our approach is
to apply a theory developed by Katok, and Katok and Mendoza
concerning the existence of hyperbolic horseshoes in the presence
of hyperbolic ergodic invariant probability measures. We show that
under the assumption of mild nonuniform hyperbolicity the
topological pressure of a   H\"older continuous potential is
entirely determined by the values of the potential on the saddle
points. Moreover, it is sufficient to consider only saddle points
with ``large" absolute value of the Lyapunov exponents. Staying in
this setting and by using saddle points, we propose  a  pressure
for non-continuous potentials and derive several variational
inequalities for it. Finally, we establish relations between the
attraction properties of the system and the dimension of the invariant set.

\subsection{Statement of the  results}
We now describe our results in more detail.
Let $M$ be a smooth Riemannian manifold and  let $\Lambda\subset M$ be a
compact
locally maximal invariant set of a $C^2$-diffeomorphism $f:M\to M$. Given
$\varphi\in C(\Lambda,\bR)$ we denote by $P_{\rm top}(\varphi)$ the
topological pressure of the potential $\varphi$, see Section~\ref{sec:2} for
the definition and details. Let $C^f(\Lambda,\bR)$ be defined as in
Section~\ref{sec:3}.
Roughly speaking a potential $\varphi$ belongs to $C^f(\Lambda,\bR)$ if
$P_{\rm top}(\varphi)$ can be approximated by  free energies of measures whose
absolute values of the Lyapunov exponents are uniformly bounded away from zero,
and $\varphi$ has no equilibrium state with zero entropy. We denote the
corresponding uniform bound by  $\delta(\varphi)$. Again we refer to
Section~\ref{sec:3} for the details. Given $n\in\bN$ let $\Fix(f^n)$ denote
the set of fixed points of $f^n$, and let $\SFix(f^n)\subset \Fix(f^n)$ denote
the saddle points in $\Fix(f^n)$. Moreover, for $0<\alpha$, $0<c\leq 1$, let
$\SFix(f^n,\alpha,c)$ be defined as in~\eqref{eqsaddle}. Roughly
speaking, $\SFix(f^n,\alpha,c)$ are those saddle  points for which
the infimum norm of the derivative restricted to the
stable/unstable spaces grows uniformly at an exponential rate
at least $\alpha$. Therefore, if $x\in \SFix(f^n,\alpha,c)$ then all
Lyapunov exponents of $x$ have absolute value greater than or equal to $\alpha$.
Let $\SPer(f)$ denote the set of all saddle points of $f$ in
$\Lambda$. The sets $\SFix(f^n,\alpha,c)$ provide a natural
filtration of $\SPer(f)$; in particular,
\begin{equation*}
\SPer(f)=\bigcup_{\alpha>0}\bigcup_{c>0}\bigcup_{n=1}^\infty \SFix(f^n,\alpha,c).
\end{equation*}
For $\varphi\in C(\Lambda,\bR)$ we define
\begin{equation}\label{eqpspal}
\P_{\rm SP}(\varphi,\alpha,c) = \limsup_{n\to\infty}
                 \frac{1}{n}\log \left(\sum_{x\in \SFix(f^n,\alpha,c)} \exp S_n\varphi(x)\right),
\end{equation}
where $S_n\varphi(x)=\sum_{k=0}^{n-1}\varphi(f^k(x))$. Our main
result shows that in the case of H\"older continuous potentials in
$C^f(\Lambda,\bR)$ the topological pressure is entirely determined by the
values of $\varphi$ on the saddle points. More precisely, we have
the following result.

\begin{theorem}\label{Main}
 Let $f\colon M\to M$ be a $C^2$-diffeomorphism and let
  $\Lambda\subset M$ be a compact locally maximal $f$-invariant set. Let
  $\varphi\in C^f(\Lambda,\bR)$ be a H\"older continuous potential and let $0<\alpha<\delta(\varphi)$. Then
  \begin{equation}\label{pressp}
  P_{\rm top}(\varphi) = \lim_{c\to 0}P_{\rm
  SP}(\varphi,\alpha,c).
  \end{equation}
\end{theorem}

It follows from Ruelle's inequality that in the case of surface
diffeomorphisms Theorem~\ref{Main} holds for a more general class
of potentials, namely~\eqref{pressp} is true for all H\"older
continuous potentials admitting no equilibrium state with zero
measure-theoretic entropy. We note that the hypothesis $\varphi\in
C^f(\Lambda,\bR)$  can, in general, not be omitted in
Theorem~\ref{Main}. A simple counterexample is given in Section~\ref{sec:3}
(see Example~\ref{ex:2}). Theorem~\ref{Main} is even new for uniformly
hyperbolic sets. Indeed, Bowen proved a version of~\eqref{pressp} in the case
of hyperbolic sets by considering 
all saddle points. On the other hand, our  result shows that it is
already sufficient to consider saddle points with ``large"
Lyapunov exponents. We emphasize, however, that the main
application of Theorem~\ref{Main} is that to non-uniformly
hyperbolic systems. Theorem~\ref{Main} is a generalization of a result of 
Chung and Hirayama~\cite{ChuHir:03}. They considered the entropy of
diffeomorphisms on compact surfaces, i.e.  $\dim M=2$, $\Lambda= M$ and
$\varphi=0$, and we use some of their ideas in our approach. It should be
noted that on higher dimensional manifolds the potential $\varphi=0$ does not
always belong to  $C^f(\Lambda,\bR)$, see Example~\ref{ex:1} in
Section~\ref{sec:3}.  

We also consider a version of pressure for
\emph{non-continuous} potentials. As a particular case we study
the volume expansion of the derivative  restricted to the
expanding subbundle $E^u$ defined by
\begin{equation*}
\varphi^u(x) = -\log\lvert\det Df(x)|E^u_x\rvert.
\end{equation*}
In general, the subbundle $E^u$ can only be defined over a particular subset
of $\Lambda$,  and $\varphi^u$ does not extend to a continuous function on
$\Lambda$. Nevertheless, since $E^u_x$ exists for all saddle points $x$, we
still can define
\begin{equation*}
P_{\rm SP}(\varphi^u)=\lim_{\alpha\to 0}\lim_{c\to 0} \P_{\rm
SP}(\varphi^u,\alpha,c)
\end{equation*}
We call $P_{\rm SP}(\varphi^u)$ the volume pressure of $f$. Note
that if $\Lambda$ is a locally maximal uniformly hyperbolic set
such that $f|\Lambda$ is topologically mixing, then by a classical
result of Bowen $P_{\rm SP}(\varphi^u)$ coincides with $P_{\rm
top}(\varphi^u)$. We show in Corollary~\ref{cor3} that in the case
of general sets $\Lambda$ which contain at least one saddle point
we have,
\begin{equation}\label{einin}
P_{\rm SP}(\varphi^u) \le
\sup_{\mu}\left(h_\mu(f)+\int_\Lambda\varphi^u d\mu\right),
\end{equation}
where the supremum is taken over all $f$-invariant probability
measures $\mu$ of saddle type. We note that even in the case of
uniformly hyperbolic sets (assuming that $f|\Lambda$ is not
topologically mixing) inequality~\eqref{einin} can be strict.

The volume pressure can be applied to characterize certain attraction
properties on $\Lambda$.  Young showed in~\cite{You:90} that
\begin{equation}\label{sais}
\sup_{\mu}\left(h_\mu(f)+\int_\Lambda\varphi^u d\mu\right) \le
\underline{E}(V) \le\overline{E}(V)\le 0,
\end{equation}
where the supremum is taken over all  ergodic invariant measures
$\mu$, and  $\underline{E}(V)$ ($\overline{E}(V)$) denotes the
lower (upper) escape rate from a neighborhood $V$ of $\Lambda$
(see~\eqref{defvol1} for precise  definition).

In the case of a locally maximal topologically mixing hyperbolic
set $\Lambda$ Bowen~\cite{B} showed (also using joint results with
Ruelle~\cite{BR}) that $P_{\rm SP}(\varphi^u)= \overline{E}(V)$;
in particular, all the inequalities in~\eqref{einin} and~\eqref{sais} are
identities and the corresponding quantity 
coincides with $P_{\rm top}(\varphi^u)$. Moreover, he proved that
$\Lambda$ is an attractor if and only if $P_{\rm
top}(\varphi^u)=0$. The latter result has recently been extended
in~\cite{ShaWol:05} by adding the  a priori weaker but still
equivalent condition $\dim_{\rm H} W^s(\Lambda)=\dim M$. Here
$\dim_{\rm H} W^s(\Lambda)$ denotes the Hausdorff dimension of the
stable set of $\Lambda$. Using ideas of~\cite{ShaWol:05} we show
for a general set $\Lambda$  that if $\overline{E}(V)<0$ then the
upper box dimension $\overline{\dim}_{\rm B}\Lambda$ of $\Lambda$
is strictly smaller than $\dim M$. More precisely, we derive in
Theorem~\ref{thSW} an upper bound for $\overline{\dim}_{\rm
B}\Lambda$ in terms of $\overline{E}(V)$ and the maximal
asymptotic exponential expansion rate of $f$ on $\Lambda$. This
bound is strictly smaller than $\dim M$ provided that
$\overline{E}(V)<0$. 

We now briefly describe the content of the paper. In
Section~\ref{sec:2} we review several concepts and results from
smooth ergodic theory and introduce various notions of  pressure.
Section~\ref{sec:3} is devoted to the statements in our main
result Theorem~\ref{Main}. In Section~\ref{secnoncon} we introduce
a version of the  pressure  for non-continuous potentials. We
study this pressure in the particular case of the potential
$\varphi^u=-\log\lvert\det Df|E^u\rvert$ and derive several
variational inequalities for it. Finally, in Section~\ref{sec:5}
we discuss relations between the escape rates and the dimension of
$\Lambda$.

\section{Preliminaries}\label{sec:2}
\subsection{Notions from smooth ergodic theory}\label{sec:2.1}

Let $M$ be a smooth Riemannian manifold and let
$f\colon M\to M$ be a $C^2$-diffeomorphism. We consider a compact locally
maximal $f$-invariant set $\Lambda\subset M$. Here \emph{locally
maximal} means that there exists an open neighborhood $U\subset M$
of $\Lambda$ such that $\Lambda=\bigcap_{n\in\bZ} f^n(U)$. Note
that in particular, if $M$ is compact, the case $\Lambda=M$ fits
within this setup. To avoid trivialities we will always assume
that $h_{{\rm top}}(f|\Lambda)>0$, where $h_{{\rm top}}$ denotes the topological entropy of the map. This rules out the case that $\Lambda$ is
only a periodic orbit. Given $x\in \Lambda$ and $v\in T_xM\backslash\{0\}$, we
define the \emph{forward Lyapunov exponent} of $v$ at $x$ (with respect to
$f$) by
\begin{equation}\label{deflya}
\lambda^+(x,v)\eqdef\limsup_{n\to\infty}\frac{1}{n}\log\lVert Df^n(x)(v)\rVert.
\end{equation}
Analogously, we define a Lyapunov exponent $\lambda^-(x,v)$ for negative time
by replacing the map $f$ in~\eqref{deflya} with $f^{-1}$, which is called the
\emph{backward Lyapunov exponent} of $v$ at $x$ (with respect to
$f$). 
Let us now assume that $x\in \Lambda$ is a \emph{Lyapunov regular} point of $f$ 
(see~\cite{BarPes:05} for the definition and details on Lyapunov regularity). 
Then there exist a positive integer $s(x)\leq \dim M$ and  a
$Df$-invariant splitting
\[
T_xM = \bigoplus_{i=1}^{s(x)}E^i_x
\]
such that for all $i=1,\ldots,s(x)$ 
and $v\in E^i_x\backslash \{0\}$ we have
\[
\lim_{n\to\pm\infty}\frac{1}{n}\log\lVert Df^n(x)(v)\rVert
 =\lambda^+(x,v) = -\lambda^-(x,v)\eqdef\lambda_i(x),
\]
with uniform convergence on $\{v\in E^i_x\colon\lVert v\rVert=1\}$.
We will count the values of the Lyapunov exponents
$\lambda_i(x)$ with their multiplicities, i.e. we consider the numbers
\[
\lambda_1(x)\le\cdots\le\lambda_{\dim M}(x).
\]

Let $\cM$ denote the set of all Borel $f$-invariant probability measures on
$\Lambda$ endowed with weak$*$ topology. This makes $\cM$
to a compact convex space. Moreover, let $\cM_{\rm E}\subset \cM$ be the subset
of ergodic measures.
By Oseledec's theorem, given $\mu\in \cM$ the set of Lyapunov regular points has full measure and
$\lambda_i(\cdot)$ is $\mu$-measurable. We denote by
\begin{equation}\label{deflyame}
\lambda_i(\mu)\eqdef\int\lambda_i(x) d\mu(x).
\end{equation}
the Lyapunov exponents of the measure $\mu$.
Note that if $\mu\in \cM_{\rm E}$ then
$\lambda_i(.)$ is constant $\mu$-a.e. and therefore, the corresponding value
coincides with $\lambda_i(\mu)$.
We say that $\mu\in\cM$ is a \emph{hyperbolic measure} if $\mu$ has non-zero
Lyapunov exponents. Set
\[
\chi(\mu)\eqdef\min_{i=1,\ldots,\dim M}\lvert\lambda_i(\mu)\rvert.
\]
In particular, if there is $1\le l=l(\mu)<\dim M$ such that
\begin{equation*}
\lambda_l(\mu)<0<\lambda_{l+1}(\mu),
\end{equation*}
we say that $\mu$ is of \emph{saddle type}. It follows from Ruelle's
inequality that for a surface diffeomorphism every hyperbolic measure $\mu$
with positive measure-theoretic entropy $h_\mu(f)$ is of saddle type. We denote
by $\Fix(f)$ the set of fixed points of $f$. Moreover, we
denote by $\Per(f)=\bigcup_n \Fix(f^n)$ the set of periodic points
of $f$. For $x\in \Fix(f^n)$ we have that $\lambda_i(x)=\frac{1}{n}\log
|\delta_i|$, where $\delta_i$ are the eigenvalues of $Df^n(x)$.
We call a periodic point $x$ a \emph{saddle point} if there is $1\le l=l(x)<\dim M$
with $\lambda_l(x)< 0 <\lambda_{l+1}(x)$. Let $\SFix(f^n)$ denote the fixed points of $f^n$ which are saddle
points. Hence,
$\SPer(f)=\bigcup_n \SFix(f^n)$ is the set of all saddle points.

We say that a compact $f$-invariant set $K\subset M$ is a \emph{hyperbolic set}
if there exists a continuous $Df$-invariant splitting of the
tangent bundle $T_K M = E^s\oplus E^u$ and constants  $c>0$
and $\lambda\in (0,1)$ such that
\begin{equation}\label{ne}
\begin{split}
\lVert Df^k(x)(v)\rVert \le c\lambda^k \lVert v\rVert \text{ for
all }v\in E^s _x,\\
\lVert Df^{-k}(x)(v)\rVert \le c\lambda^k \lVert v\rVert \text{
for all }v\in E^u_x
\end{split}
\end{equation}
for all $x\in K$ and all $k\in\bN$.
For convenience we sometimes also refer to relative compact sets satisfying
\eqref{ne} as hyperbolic sets.

For $x\in \SFix(f^n)$ we define $E^s_x$ ($E^u_x$) to be the direct sum of
the eigenspaces corresponding to eigenvalues of $Df^n(x)$ with norm smaller
than 1 (larger than 1).
It is easy to see that there exists  $0<c\leq 1$, $c=c(x)$ such that for all
integers $k\ge 0$ and $0\le i\le n-1$
\begin{equation}\label{eqhi}
\begin{split}
c e^{k\lambda_{l+1}(x)} &\le \lVert Df^k(f^i(x)) (v) \rVert\\
c e^{-k\lambda_l(x)} &\le \lVert Df^{-k}(f^i(x)) (w)\rVert
\end{split}
\end{equation}
whenever $v\in E^u_{f^i(x)}$ with $\lVert v\rVert=1$ and $w\in E^s_{f^i(x)}$
with $\lVert w\rVert=1$.
For $0<\alpha$, $0<c\leq 1$, and $n\in\bN$ we set
\begin{multline}\label{eqsaddle}
\SFix(f^n,\alpha,c) \eqdef\{ x\in\SFix(f^n)\colon
\lVert(Df^{-k}(f^i(x))|E^s)^{-1}\rVert^{-1}\ge c e^{k\alpha},\\
\lVert(Df^k(f^i(x))|E^u)^{-1}\rVert^{-1} \ge c e^{k\alpha}
 \text{ for all } k\ge 1 \text{ and } 0\le i\le n-1 \}.
\end{multline}
Thus, if  $\alpha\geq \alpha',\,c\geq c'$, then
\begin{equation}\label{ni}
\SFix(f^n,\alpha,c) \subset \SFix(f^n,\alpha',c')
\end{equation}
and
\begin{equation*}
\SPer(f)=\bigcup_{\alpha>0}\bigcup_{c>0}\bigcup_{n=1}^\infty
\SFix(f^n,\alpha,c).
\end{equation*}

\begin{lemma}\label{lemma1}
  Let  $\alpha$, $c>0$ be fixed. Then $x\mapsto E^{s/u}_x$ are continuous
   maps on $L= \bigcup_{n=1}^\infty \SFix(f^n,\alpha,c)$ which extend
  continuously to the closure of $L$.
\end{lemma}

\begin{proof}
  It suffices to notice that for $y\in L$ we have
  \[
  \lVert Df^k(y)|E^s\rVert\lVert Df^{-k}(f^k(y))|E^u\rVert \le
  \frac{1}{c^2}e^{-2k\alpha} \le \frac{1}{2}
  \]
  whenever $k\ge \frac{1}{2\alpha}\log \frac{2}{c^2}$. This means that the
  splitting of the tangent space $T_LM=E^s\oplus E^u$ is $k$-dominated for any
  such $k$ (see~\cite{BonDiaPuj:03} for details on dominated splittings). By
  $f$-invariance of $L$ the statement follows from~\cite[Lemma 1.4]{BonDiaPuj:03}.
\end{proof}

We now present some results concerning non-uniformly hyperbolic
systems developed by Katok and Mendoza, see~\cite{Kat:80} in the
case of surface diffeomorphisms and~\cite{BarPes:05} in the general
case. Let $\mu\in \cM_E$ be a hyperbolic measure with positive
measure-theoretic entropy. Then
\begin{equation}\label{eqmesent}
h_\mu(f)\le \limsup_{n\to\infty}\frac{1}{n}\log \card\Fix(f^n).
\end{equation}
Moreover, for all $\varepsilon>0$ we have
\begin{equation*}
0<h_\mu(f) \le \limsup_{n\to\infty}\frac{1}{n}\log
\card\left\{x\in\SFix(f^n)\colon\chi(x)\ge \chi(\mu) -\varepsilon
\right \},
\end{equation*}
where $\chi(x)=\min_{j=1,\ldots,\dim M}\lvert\lambda_j(x)\rvert$.
Furthermore, there exist positive constants $\alpha_0$ and $c_0$
such that $\SFix(f^n,\alpha,c)\not=\emptyset$ for all
$\alpha\le\alpha_0$, $c\le c_0$, and infinitely many $n\in\bN$.

\subsection{Various  pressures}
Next,  we introduce a version of topological pressure which is
entirely determined by the values of the potential on the saddle
points.

Let us first recall the classical topological pressure. Let
$(\Lambda,d)$ be a compact metric space and let $f\colon
\Lambda\to \Lambda$ be a continuous map. For $n \in {\mathbb N}$ we
define a new metric $ d_n $ on $ \Lambda$ by $
d_n(x,y)=\max_{k=0,\ldots ,n-1} d(f^k(x),f^k(y))$. A set of points
$\{ x_i\colon i\in I \}\subset \Lambda$ is called
\emph{$(n,\varepsilon)$-separated} (with respect to $f$) if
$d_n(x_i,x_j)> \varepsilon$ holds for all $x_i,x_j$ with $x_i \ne
x_j$. Fix for all $\varepsilon>0$ and all $n\in\bN$ a maximal
(with respect to the inclusion) $(n,\varepsilon)$-separated set
$F_n(\epsilon)$. The \emph{topological pressure} (with respect to $f|\Lambda$) is a mapping
$ P_{\rm top}(f|\Lambda,.)\colon C(\Lambda,\bR)\to \bR$  defined by
\begin{equation}\label{defdru}
  P_{\rm top}(f|\Lambda,\varphi) \eqdef \lim_{\varepsilon \to 0}
            \limsup_{n\to \infty}
            \frac{1}{n} \log \left(\sum_{x\in F_n(\epsilon)}
            \exp S_n\varphi(x) \right),
\end{equation}
where
\begin{equation}\label{eqsn}
S_n\varphi(x)\eqdef\sum_{k=0}^{n-1}\varphi(f^k(x)).
\end{equation}
The \emph{topological entropy} of $f$ on $\Lambda$ is defined by
$h_{\rm top}(f|\Lambda)=P_{\rm top}(f|\Lambda,0)$. We simply write $P_{\rm
  top}(\varphi)$ and $h_{\rm top}(f)$ if there is no confusion about $f$ and
$\Lambda$.
Note that the definition of $P_{\rm top}(\varphi)$ does not depend on the choice of the
sets $F_n(\epsilon)$ (see~\cite{Wal:81}).
The topological pressure satisfies the
following variational principle:
\begin{equation}\label{eqvarpri}
P_{\rm top}(\varphi)= \sup_{\nu\in \cM} \left(h_\nu(f)+\int_\Lambda \varphi
\,d\nu\right).
\end{equation}
Furthermore, the supremum in~\eqref{eqvarpri} can be replaced by
the supremum taken only over all $\nu\in\cM_{\rm E}$. We now
introduce a \emph{pressure} which is entirely defined  by the
values of $\varphi$ on the saddle points. Let $\varphi\in
C(\Lambda, \bR)$ and let $0<\alpha$, $0<c\leq 1$.
Define
\begin{equation*}
 Q_{\rm SP}(\varphi,\alpha,c,n) \eqdef
\sum_{x\in \SFix(f^n,\alpha,c)} \exp S_n\varphi(x)
\end{equation*}
if $\SFix(f^n,\alpha,c)\ne\emptyset$ and
\begin{equation*}
Q_{\rm SP}(\varphi,\alpha,c,n)\eqdef
\exp\left(n\min_{x\in\Lambda} \varphi(x)\right)
\end{equation*}
otherwise. Furthermore, we define
\begin{equation*}
P_{\rm SP}(\varphi,\alpha,c) \eqdef \limsup_{n\to\infty}
                 \frac{1}{n}\log Q_{\rm SP}(\varphi,\alpha,c,n).
\end{equation*}
It follows from the  definition that if
$\SFix(f^n,\alpha,c)\not=\emptyset$ for some $n\in\bN$
then this is true already for infinitely many $n\in\bN$. Therefore,
in the case when $\SFix(f^n,\alpha,c)\not=\emptyset$ for
some $n\in\bN$ then $P_{\rm SP}(\varphi,\alpha,c)$ is
entirely determined by the values of $\varphi$ on
$\bigcup_{n\in\bN}\SFix(f^n,\alpha,c)$.

The following classical result shows that in the case of
hyperbolic sets the topological pressure of a H\"older continuous
potential is entirely determined by its values on the saddle
orbits. We only sketch its proof and refer to~\cite{KatHas:95} for
full details.

\begin{proposition}\label{ha}
 Let $f\colon M\to M$ be a $C^2$-diffeomorphism and let $K\subset M$ be
 a compact hyperbolic  set of $f$. Let $\varphi\in C(K,\bR)$ be a
 H\"older continuous potential.  Then
 \begin{equation}\label{eqdruck}
 \limsup_{n\to\infty}\frac{1}{n}\log\left(\sum_{x\in\Fix(f^n)\cap K}\exp
 S_n\varphi(x)\right) \le P_{\rm top}(f|K,\varphi).
 \end{equation}
 Furthermore, if $f|K$ satisfies the specification property then we have
 equality in~\eqref{eqdruck}, and the limit superior is in fact a limit.
\end{proposition}

\begin{proof}
Since $K$ is a hyperbolic set the map $f|K$ is expansive. If
$\delta$ is the expansivity constant then for every $n\in\bN$ and
every $0<\varepsilon\le \delta$ the set $\Fix(f^n)$ is
$(n,\varepsilon)$-separated. Thus, the inequality~\eqref{eqdruck} follows from the fact that
the definition~\eqref{defdru} can be replaced by the supremum taken over all
$(n, \epsilon)$-separated sets (see~\cite{Wal:81}).
The equality in~\eqref{eqdruck} is a direct consequence of~\cite[Proposition
20.3.3]{KatHas:95}.
\end{proof}

{\it Remark. } It follows from Proposition~\ref{ha} and the Specification Theorem (see~\cite[Theorem
18.3.9]{KatHas:95}) that if $K$ is a locally maximal hyperbolic set such that $f|K$ is topologically mixing then~\eqref{eqdruck} is an identity.

\section{Saddle points and topological pressure}\label{sec:3}
In this section we study possible extensions of Proposition~\ref{ha}. In
particular, we consider  general locally maximal 
invariant sets $\Lambda$ without requiring hyperbolicity. Let $M$
be a smooth Riemannian manifold and let $f\colon M\to M$ be a
$C^2$-diffeomorphism. Let $\Lambda\subset M$ be a compact
 locally maximal $f$-invariant set. Let $\varphi\in
C(\Lambda,\bR)$. First, we investigate as to how periodic points
with small absolute value of the Lyapunov exponents contribute to
the pressure $P_{\rm top}(\varphi)$. For $0<\alpha<\beta$ and
$0<c\le 1$ we define
\[
\begin{split}
\SFix&(f^n,[\alpha,\beta],c) =\{ x\in\SFix(f^n)\colon \\
&ce^{k\alpha}\le
\lVert Df^k(f^i(x))(v)\rVert,
\lVert Df^{-k}(f^i(x))(w)\rVert \le
c^{-1} e^{k\beta}\\
&\text{ for all }k\ge 0\text{ and }0\le i\le n-1 \text{ and all }v\in
E^u_x,w\in E^s_x\}.
\end{split}
\]
It follows immediately from~\eqref{eqhi}  that for every $x\in
\SFix(f^n,[\alpha,\beta],c)$ we have $\alpha \le
\lvert\lambda_i(x)\rvert\le \beta$  for all $i=1,\ldots,\dim M$.
Furthermore,
\begin{equation}\label{he}
\SFix(f^n,\alpha,c) =\bigcup_{\beta>0}\SFix(f^n,[\alpha,\beta],c)
                     =\SFix(f^n,[\alpha,\beta_0],c),
\end{equation}
where
\begin{equation}\label{beta}
  \beta_0\eqdef \max\{\log\lVert Df(x)\rVert, \log\lVert Df^{-1}(x)\rVert\colon
  x\in \Lambda\}.
\end{equation}
Define
\begin{equation*}
 Q_{\rm SP}(\varphi,[\alpha,\beta],c,n) =
\sum_{x\in \SFix(f^n,[\alpha,\beta],c)} \exp S_n\varphi(x)
\end{equation*}
if $\SFix(f^n,[\alpha,\beta],c)\ne\emptyset$ and
\begin{equation}\label{part}
Q_{\rm SP}(\varphi,[\alpha,\beta],c,n)=
\exp\left(n\min_{x\in\Lambda} \varphi(x)\right)
\end{equation}
 otherwise.
Furthermore, we define
\begin{equation*}
P_{\rm SP}(\varphi,[\alpha,\beta],c) =
\limsup_{n\to\infty}\frac{1}{n}\log Q_{\rm
SP}(\varphi,[\alpha,\beta],c,n) .
\end{equation*}
It follows from the  definition that if
$\SFix(f^n,[\alpha,\beta],c)\not=\emptyset$ for some $n\in\bN$
then this is already true  for infinitely many $n\in\bN$. Therefore,
in the case when $\SFix(f^n,[\alpha,\beta],c)\not=\emptyset$ for
some $n\in\bN$ then $P_{\rm SP}(\varphi,[\alpha,\beta],c)$ is
entirely determined by the values of $\varphi$ on
$\bigcup_{n\in\bN}\SFix(f^n,[\alpha,\beta],c)$. We have the
following:

\begin{theorem}\label{theo2}
  Let $f\colon M\to M$ be a $C^2$-diffeomorphism and let $\Lambda\subset
  M$ be a compact locally maximal $f$-invariant  set. Let
  $0<\alpha<\beta$ and $0<c\leq 1$ such that
  $\SFix(f^n,[\alpha,\beta],c)\not=\emptyset$ for some $n\in\bN$.
  Let $\varphi\in C(\Lambda,\bR)$ be a H\"older continuous potential. Then
  \begin{equation}\label{le}
    P_{\rm SP}(\varphi,[\alpha,\beta],c) \le
    \sup_\nu \left\{h_\nu(f)+\int_\Lambda \varphi d\nu \right\}
    \leq P_{\rm top}(\varphi),
  \end{equation}
  where the supremum is taken over all $\nu\in \cM_{\rm E}$
  with $\alpha\le\lvert\lambda_i(\nu)\rvert\le\beta$ for all $i=1,\ldots,\dim
  M$.
\end{theorem}

\begin{proof}
Let $0<\alpha<\beta$ and $0<c\leq 1$ such that
$\SFix(f^n,[\alpha,\beta],c)\not=\emptyset$ for some $n\in\bN$.
In particular, the supremum in~\eqref{le} is not taken over the empty set.
The right hand side  inequality in~\eqref{le} is
a consequence of the variational principle.
In order to prove the left hand side inequality set
  \[
  K=K_{\alpha,\beta,c} \eqdef \overline{\bigcup_{n=1}^\infty
  \SFix(f^n,[\alpha,\beta],c)} .
  \]
  The subspaces of the $Df$-invariant splitting $T_x M=E^s_x\oplus E^u_x$ vary
  continuously on the set $\bigcup_{n=1}^\infty \SFix(f^n,[\alpha,\beta],c)$
  and by Lemma~\ref{lemma1} they can be extended continuously to $K$.
  It follows that $K$ is a  hyperbolic set for $f$.
  Furthermore, for every $n\ge 1$ with $\SFix(f^n,[\alpha,\beta],c)\not=\emptyset$ we have,
  \begin{equation}\label{eqsi}
  \Fix(f^n)\cap K = \SFix(f^n,[\alpha,\beta],c).
  \end{equation}
  Therefore, Proposition~\ref{ha} implies
\begin{equation}\label{ghj}
P_{\rm SP}(\varphi,[\alpha,\beta],c)\leq P_{\rm top}(f|K,\varphi).
\end{equation}
It follows from  the variational principle that for every
$\varepsilon>0$ there is a $\mu\in \cM_E$ which is supported
  in
  $K$ such that
  \begin{equation}\label{hjk}
  P_{\rm top}(f|K,\varphi) -\varepsilon
  \le h_\mu(f) + \int_K\varphi d\mu
  \le P_{\rm top}(f|K,\varphi).
  \end{equation}
Since $\mu$ is ergodic we have that $\lambda_i(x)=\lambda_i(\mu)$
for $\mu$-almost every $x\in K$.  It now follows from
 the continuity of $x\mapsto Df(x)$, the continuity of the
 extended splitting $T_K M= E^s\oplus E^u$ and the definition of
$\SFix(f^n,[\alpha,\beta],c)$ that $\alpha\leq |\lambda_i(x)|\leq
\beta$ for all $x\in K$ and all $i=1,\ldots,\dim M$. We conclude
that $\alpha\le\lvert\lambda_i(\mu)\rvert\le\beta$ for all
$i=1,\ldots,\dim M$. Therefore, the left hand side inequality in
\eqref{le} follows from~\eqref{ghj} and~\eqref{hjk}. This
completes the proof.
\end{proof}

\begin{proposition}\label{li}
  Let $f\colon M\to M$ be a $C^2$-diffeomorphism and let $\Lambda\subset M$ be
  a compact locally maximal $f$-invariant set. Let $\varphi\in
  C(\Lambda,\bR)$ be a H\"older continuous potential.
  Then for all $\mu\in \cM_{\rm E}$ with $h_\mu(f)>0$ and $\chi(\mu)>0$ and for
  all $0<\alpha<\chi(\mu)$ we have
  \begin{equation*}
  h_\mu(f)+\int_\Lambda\varphi d\mu \le
  \lim_{c\to 0} P_{\rm SP}(\varphi,\alpha,c).
  \end{equation*}
\end{proposition}

\begin{proof}
Consider $\mu\in \cM_{\rm E}$ with $h_\mu(f)>0$ and $\chi(\mu)>0$,
and let $0<\alpha<\chi(\mu)$. It follows from the work of Katok
and Mendoza (see e.g.~\cite[Chapter S.5]{KatHas:95} for the case
of surface diffeomorphisms and~\cite{BarPes:05} for the general
case) that there exists a sequence $(\mu_n)_n$ of measures
$\mu_n\in \cM_{\rm E}$ supported on hyperbolic horseshoes
$K_n\subset M$ (see~\cite{KatHas:95} for the definition) such that
$\mu_n\to\mu$ in the weak$\ast$ topology, $h_{\mu_n}(f)\to
h_\mu(f)$, and $\chi(\mu_n)\to\chi(\mu)$. In particular, for each
$n\in\bN$ there exist $m,s\in\bN$ such that $f^m| K_n$ is
conjugate to the full shift in $s$ symbols. Since $\Lambda$ is a
compact locally maximal $f$-invariant set we can conclude that
$K_n\subset \Lambda$ for all $n\in \bN$. It follows that  for
every $0<\varepsilon<\chi(\mu)-\alpha$ there is a number
$n=n(\varepsilon)\ge 1$ such that
\begin{equation}\label{nag}
  h_\mu(f)-\varepsilon < h_{\mu_n}(f) \text{ and }
  \int_\Lambda\varphi d\mu-\varepsilon < \int_\Lambda\varphi d\mu_n.
\end{equation}
Moreover, there exists a number $c_0=c_0(n)$ with $0<c_0(n)\le1$ such that for
every periodic point $x\in K_n$ and every $k\in\bN$ we have
\[
\begin{split}
  \lVert Df^k(x)(v)\rVert &\ge  c_0e^{k(\chi(\mu)-\varepsilon)}\lVert v\rVert
  \text{ for every }v\in E^u_x\\
  \lVert Df^{-k}(x)(w)\rVert &\le c_0e^{k(-\chi(\mu)+\varepsilon)}\lVert
  w\rVert
  \text{ for every }w\in E^s_x.
\end{split}
\]
This implies that
\begin{equation}\label{kir}
  \Fix(f^k)\cap K_n \subset \SFix(f^k,\alpha,c_0)
\end{equation}
for every $k\in\bN$. It follows from~\eqref{nag} and the variational
principle~\eqref{eqvarpri} that
\begin{equation*}
  h_\mu(f)+\int_\Lambda\varphi d\mu -2\varepsilon
  < h_{\mu_n}(f)+\int_{K_n}\varphi d\mu_n
  \le P_{\rm top}(f|K_n,\varphi).
\end{equation*}
Let $m, s\in\bN$ such that
$f^m|K_n$ is topologically conjugate to the full shift in $s$
symbols.
Since $mP_{\rm top}(f|K_n,\varphi) = P_{\rm top}(f^m|K_n,S_m\varphi)$
(see~\cite[Theorem 9.8]{Wal:81}), we can conclude that
\begin{equation*}
  h_\mu(f)+\int_\Lambda\varphi d\mu -2\varepsilon
  \le \frac{1}{m}P_{\rm top}(f^m|K_n,S_m\varphi).
\end{equation*}
Recall that $S_m\varphi(x)=\sum_{i=0}^{m-1}\varphi(f^i(x))$.
It now follows from Proposition~\ref{ha} and an elementary calculation that
\begin{equation}\label{eqrep}
  \begin{split}
    h_\mu(f)&+\int_\Lambda\varphi d\mu -2\varepsilon\\
    &\le \frac{1}{m}\lim_{k\to\infty}\frac{1}{k}\log\left(\sum_{x\in\Fix(f^{mk})\cap K_n}
      \exp \left(\sum_{i=0}^{k-1}S_m\varphi(f^{im}(x))\right)\right)\\
    &=\lim_{k\to\infty}\frac{1}{mk}\log\left(\sum_{x\in\Fix(f^{mk})\cap K_n}
      \exp S_{mk}\varphi(x)\right)\\
    &\le \lim_{k\to\infty}\frac{1}{k}\log\left(\sum_{x\in\Fix(f^{k})\cap K_n}
      \exp S_{k}\varphi(x)\right).
  \end{split}
\end{equation}
Combining~\eqref{kir} and~\eqref{eqrep} yields
\begin{equation*}
  h_\mu(f)+\int_\Lambda\varphi d\mu -2\varepsilon
  \le \limsup_{k\to\infty}\frac{1}{k}\log
  \sum_{x\in\SFix(f^k,\alpha,c_0)}\exp S_k\varphi(x).
\end{equation*}
Recall that by~\eqref{ni} the  map $c\mapsto P_{\rm SP}(\varphi,\alpha,c)$ is
non-decreasing as $c\to 0^+$.
Since $\varepsilon>0$ was chosen arbitrarily  the claimed statement follows.
\end{proof}

\noindent
{\it Remarks. }\\
(i) We note that the hyperbolic horseshoes $K_n$ in the proof of Proposition~\ref{li} are in general not locally maximal $f$-invariant sets.\\
(ii) It follows from Ruelle's inequality that if $M$ is a surface
then $h_\mu(f)>0$ implies $\chi(\mu)>0$ and therefore, Proposition
\ref{li} holds for all measures with positive measure-theoretic
entropy.
\\[0.3cm]
\noindent
We now introduce a natural class of potentials. For $\varphi\in
C(\Lambda,\bR)$ set
\begin{equation}\label{defal}
\alpha(\varphi)=
P_{\rm top}(\varphi)-\sup_{\nu\in\cM}\int_\Lambda\varphi
d\nu.
\end{equation}
We say that a potential $\varphi$ belongs to $C^f(\Lambda,\bR)$ if
\begin{enumerate}
\item [(a)] $\alpha(\varphi)>0$; \item [(b)] there exist $0<\delta(\varphi)<\alpha(\varphi)$ and a
sequence $(\mu_n)_n\subset \cM_{\rm E}$ 
   such that
  $\chi(\mu_n)>\delta(\varphi)$ for every $n\in\bN$ and
  $h_{\mu_n}(f)+\int_\Lambda \varphi d\mu_n \to P_{\rm top}(\varphi)$ as
  $n\to\infty$.
\end{enumerate}
\ni {\it Remarks.}\\
(i) Note that $\alpha(\varphi)\geq 0$, and $\alpha(\varphi)>0$ if and only if
$\varphi$ has no equilibrium state with zero entropy. \\
(ii) It follows from Ruelle's inequality and the variational
principle that if $M$ is a surface then
property (b)  follows from property (a).\\
(iii) We note that if $\varphi$ has a hyperbolic equilibrium
measure $\mu_\varphi$ then we can simply choose  the constant sequence
$\mu_n=\mu_\varphi$ in (b).\\
(iv) It is easy to see that on higher dimensional manifolds the potential $\varphi=0$ does
not always belong to $C^f(\Lambda,\bR)$. A simple counterexample is given below.

\begin{example}\label{ex:1}{\rm
    Let $f=(f_1,f_2):S^1\times T^2 \to S^1\times T^2$, where $f_1$ is a rotation
    on the circle and $f_2$ is an Anosov diffeomorphism on the 2-torus. Set
    $\Lambda=S^1\times T^2$ and $\varphi=0$. Then $\alpha(\varphi)=h_{\rm
      top}(f)=h_{\rm top}(f_2)>0$; nevertheless property (b) is not satisfied
    since $f$ does not have any hyperbolic measure.
}\end{example}

We now present the proof of our main result Theorem~\ref{Main}  which was
stated in the introduction.
\begin{proof}[Proof of Theorem 1]
Let $0<\alpha<\delta(\varphi)$ and
$0<\epsilon<\delta(\varphi)-\alpha$. Let $\beta_0$ be defined as
in~\eqref{beta}. It follows from Theorem~\ref{theo2},~\eqref{he}
and~\eqref{part} that
 \begin{equation*}
 P_{\rm SP}(\varphi,\alpha,c) = P_{\rm SP}(\varphi,[\alpha,\beta_0],c) \le
 P_{\rm top}(\varphi)
 \end{equation*}
 for all $0<c\le1$.
 It remains to prove that
 \begin{equation*}
 P_{\rm top}(\varphi)\le \lim_{c\to 0}P_{\rm SP}(\varphi,\alpha,c).
 \end{equation*}
 Since $\varphi\in C^f(\Lambda,\bR)$, there exist $n\in \bN$ and
 $\mu_n\in\cM_{\rm E}$ with $\chi(\mu_n)>\delta(\varphi)$ such that
 \begin{equation}\label{gu}
 P_{\rm top}(\varphi)-\varepsilon
 \le h_{\mu_n}(f) + \int_\Lambda\varphi d\mu_n.
 \end{equation}
 It follows from~\eqref{defal} that
\begin{equation*}
 \alpha<P_{\rm top}(\varphi) -\sup_{\nu\in\cM}\int_\Lambda\varphi d\nu -\varepsilon
 \le h_{\mu_n}(f) + \int_\Lambda\varphi d\mu_n
     -\sup_{\nu\in\cM}\int_\Lambda\varphi d\nu
 \le h_{\mu_n}(f).
 \end{equation*}
 Therefore,  Proposition~\ref{li} implies
 \begin{equation}\label{eqwer}
   h_{\mu_n}(f)+\int_\Lambda\varphi d\mu_n
   \le \lim_{c\to 0}P_{\rm SP}(\varphi,\alpha,c).
 \end{equation}
 Since $\varepsilon$ can chosen arbitrary small,~\eqref{gu} and~\eqref{eqwer}
 imply~\eqref{pressp}.  Finally, if $M$ is a  surface then by Ruelle's
 inequality property (b) in the definition of $C^f(\Lambda,\bR)$ holds for all $0<\delta(\varphi)<\alpha(\varphi)$. Thus,
\eqref{pressp} holds for all $0<\alpha<\alpha(\varphi)$.
\end{proof}
\ni {\it Remarks. }\\
(i) Note that Theorem~\ref{Main} holds for every
fixed $\alpha$ satisfying $0<\alpha<\delta(\varphi)$. Therefore,
the topological pressure of $\varphi$ is entirely determined by
the values of $\varphi$ on the saddle
points of $f$.\\
(ii) As stated in the introduction Theorem 1 (and also Theorem 2) generalizes
 results of Chung and Hirayama~\cite{ChuHir:03}.  They considered the entropy (i.e. $\varphi=0$) of surface diffeomorphisms, and our results allow more general potentials in $C^f(\Lambda,\bR)$ and also diffeomorphisms on
higher dimensional manifolds. We note that we have used several ideas of~\cite{ChuHir:03} in our proofs.
For related results in the case of $C^1$-maps on the interval or the circle we refer
to~\cite{Chu}.
\\[0.3cm]
\ni
 The following example shows that the hypothesis
$\alpha(\varphi)>0$ in Theorem 1 can not be omitted.

\begin{example}\label{ex:2}{\rm 
    Let $f\colon \bR^2\to\bR^2$ be a
    $C^2$-diffeomorphism having a hyperbolic horseshoe
    $\widetilde\Lambda$ as well as an attracting fixed point $x$.
    Define $\Lambda=\widetilde\Lambda\cup\{x\}$. Moreover, let
    $\varphi\in C(\Lambda,\bR)$ with $\varphi|\widetilde\Lambda=0$ and
    $\varphi(x)> h_{\rm top} (f|\Lambda)$. Then $\alpha(\varphi)= 0$
    and $P_{\rm top}(\varphi)=\varphi(x) > \lim_{c\to 0}P_{\rm
      SP}(\varphi,\alpha,c)=h_{\rm top}(f|\Lambda)$ for every $0<\alpha<h_{\rm
      top}(f|\Lambda)$.
}\end{example}

\section{Non-continuous potentials -- the volume pressure}\label{secnoncon}

It is well-known that if $\Lambda$ is a locally maximal hyperbolic set of a
$C^2$-diffeomorphism $f\colon M\to M$ then the value of the topological
pressure of
the potential $\varphi^u=-\log|\det Df|E^u|$ has significant impact on the
geometry of $\Lambda$ as well as on the dynamics of $f$ in a neighborhood of
$\Lambda$. For example, the classical result of Bowen~\cite{B} states that
$P_{\rm top}(\varphi^u)=0$ if and only if $\Lambda$ is an attractor. In the
case of more general systems, no such continuous $Df$-invariant splitting may
exist, and therefore, one can not even define the  continuous potential
$\varphi^u$ on $\Lambda$.
To overcome this problem we introduce in this section a ``potential"
$\varphi^u$ defined only on a certain subset of $\Lambda$ and then consider
the saddle point pressure of $\varphi^u$  rather than its topological
pressure.

Let $M$ be a smooth Riemannian manifold and let $f\colon M\to M$ be a
$C^2$-diffeomorphism. Suppose that $\Lambda\subset M$ is a compact locally
maximal $f$-invariant set with $h_{\rm top} (f|\Lambda)>0$. Define
\begin{equation*}
\chi(\Lambda)\eqdef\sup_\mu \chi(\mu),
\end{equation*}
where the supremum is taken over all $\mu\in \cM_E$ supported in $\Lambda$
which have at least one positive as well as one negative Lyapunov
exponent. Since $h_{\rm top} (f|\Lambda)>0$, the variational principle for the
topological entropy and Ruelle's inequality imply that the supremum in the
definition of $\chi(\Lambda)$ is not taken over the empty set.
Note that $\chi(\Lambda)>0$ is equivalent to the existence of a  measure
$\mu\in \cM_E$ of saddle type. We denote by $\cR$ the set of  Lyapunov
regular points in $\Lambda$ (see Section~\ref{sec:2.1} and~\cite{BarPes:05} for  details). Moreover, we define
\begin{equation*}
\cR_\pm \eqdef\{x\in\cR\colon
\lambda_{1}(x)<0<\lambda_{\dim M}(x)\}
\end{equation*}
and
\begin{equation*}
\cR_{\rm H} \eqdef \{x\in\cR\colon
\lambda_{l(x)}(x)<0<\lambda_{l(x)+1}(x)\text{ for some }1\le l(x)< \dim M\}.
\end{equation*}
Clearly $\cR_{\rm H}\subset \cR_\pm$.
Let  $\mu\in \cM_E$ have at least one positive and one negative Lyapunov exponent.
It follows from the multiplicative ergodic theorem that
\begin{equation}\label{reg}
\lambda_i(x)=\lambda_i(\mu)
\end{equation}
for $\mu$-almost every $x\in\Lambda$.
Denote by $\Lambda(\mu)$ the set of points in $\cR_{\pm}$
satisfying~\eqref{reg}. Hence $\mu(\Lambda(\mu))=1$ and, in particular, if
$\mu$ is of saddle type then $\Lambda(\mu)\subset \cR_{\rm H}$. Given $x\in
\cR_\pm$ we define subspaces
\begin{equation*}
E^s_x\eqdef\bigoplus_{i\colon \lambda_i(x)<0} E^i_x,\quad
E^c_x\eqdef\bigoplus_{i\colon \lambda_i(x)=0} E^i_x\quad\text{ and }
\quad
E^u_x\eqdef\bigoplus_{i\colon \lambda_i(x)>0} E^i_x.
\end{equation*}
Hence, $x\mapsto E^{u/s/c}_x$ are Borel measurable functions on $\Lambda(\mu)$
and form a $Df$-invariant splitting of the tangent bundle
$T_{\Lambda(\mu)}M=E^s\oplus E^c \oplus E^u$. Moreover, if $\mu$ is of saddle
type then $E^c_x=\{0\}$  for $\mu$-almost every $x\in \Lambda$.
We define
\begin{equation*}
\cL \eqdef \bigcup_{\mu}\Lambda(\mu),
\end{equation*}
where the union is taken over all measures $\mu\in \cM_E$ having at least one positive and one negative Lyapunov exponent.
Analogously, we define  $\cL_H = \bigcup_{\mu}\Lambda(\mu)$, where the union is taken over all measures $\mu\in \cM_E$ of saddle type.
We define  $\varphi^u\colon \cL\to\bR$ by
\begin{equation*}
\varphi^u(x)\eqdef -\log\lvert\det Df(x)|E^u_x\rvert.
\end{equation*}
We note that in general neither  the splitting \(T_{\cL}M=E^s\oplus E^c\oplus E^u\), nor the splitting \(T_{\cL_H}M=E^s\oplus E^u\)
can  be  continuously extended to $T_\Lambda M$, and therefore there is no
continuous function \(\varphi\in C(\Lambda,\bR)\) with
$\varphi|\cL = \varphi^u$ or $\varphi|\cL_H = \varphi^u$.
Clearly, $\SPer(f) \subset \cL_H$,
and thus, $\varphi^u(x)$ is well-defined for every $x\in \SPer(f)$.
Let $0<\alpha$ and $0<c\leq 1$.
Define
\begin{equation*}
  Q_{\rm SP}(\varphi^u,\alpha,c,n) \eqdef
  \sum_{x\in \SFix(f^n,\alpha,c)} \exp S_n\varphi^u(x)
\end{equation*}
if $\SFix(f^n,\alpha,c)\ne\emptyset$ and
\begin{equation*}
  Q_{\rm SP}(\varphi^u,\alpha,c,n)\eqdef
  \exp\left(n\inf_{x\in\cL} \varphi^u(x)\right)
\end{equation*}
otherwise. We define
\begin{equation*}
  P_{\rm SP} (\varphi^u,\alpha,c) \eqdef
  \limsup_{n\to\infty}\frac{1}{n}\log Q_{\rm SP}(\varphi^u,\alpha,c,n).
\end{equation*}
It follows that if $\SFix(f^n,\alpha,c)\not=\emptyset$ for some $n\in\bN$ then
$P_{\rm SP} (\varphi^u,\alpha,c)$ is entirely determined by  the values of
$\varphi^u$ on the saddle points of $f$.
By $Df$-invariance of $E^u$ on $\cL$, we conclude that
\begin{equation*}
 \lvert \det Df^n(x)|E^u_x\rvert^{-1} = \exp S_n\varphi^u(x)
\end{equation*}
for all $x\in \cL$. Moreover, if $x\in\cL$ then
\begin{equation*}
  \lim_{n\to\infty}\frac{1}{n}\log\lvert\det Df^n(x)|E^u_x\rvert =
  \sum_{i\colon\lambda_i(x)>0} \lambda_i(x).
\end{equation*}
In particular,
\begin{equation*}
  \int_\Lambda \varphi^ud\mu =
  - \int_\Lambda \sum_{i\colon\lambda_i(x)>0} \lambda_i(x) d\mu(x)
\end{equation*}
for every measure $\mu\in\cM_{\rm E}$ of saddle type.

We now state a variational inequality for the above pressure.

\begin{theorem}\label{varprinc}
  Let \(f\colon M\to M\) be a \(C^2\)-diffeomorphism and let $\Lambda\subset
  M$ be a compact locally maximal $f$-invariant  set. Let $0<\alpha$ such that
  $\SFix(f^n,\alpha,c_0)\not=\emptyset$ for some $n\in \bN$ and some $0<c_0\leq
  1$. Then
  \begin{equation}\label{flu}
    \lim_{c\to 0}P_{\rm SP} (\varphi^u,\alpha,c)
    \le\sup_\nu \left( h_\nu(f)+\int_\Lambda \varphi^ud\nu\right),
  \end{equation}
   where the supremum is taken over all
  measures $\nu\in\cM_{\rm E}$ of saddle type with $\alpha\le\chi(\nu)$.
\end{theorem}

\begin{proof}
  Since  there is $n\in\bN$ and $0<c_0\leq 1$ such that
  $\SFix(f^n,\alpha,c_0)\not=\emptyset$, it follows that $\chi(\Lambda)>0$; in
  particular, the supremum in~\eqref{flu} is not taken over the empty set. Let
  $\beta_0$ be defined as in~\eqref{beta} and let $0<c\leq c_0$.
  Analogously as in the proof of Theorem~\ref{theo2} we are able to construct a
  compact hyperbolic  set $K=K_{\alpha,\beta_0,c}\subset \Lambda$ with
  \begin{equation}\label{fix}
    \Fix(f^n)\cap K = \SFix(f^n,\alpha,c)
  \end{equation}
  for all $n\in \bN$. Hence
  \begin{equation}\label{fix2}
    P_{\rm SP} (\varphi^u,\alpha,c) =
    \limsup_{n\to\infty}\frac{1}{n}\log\left(\sum_{x\in\Fix(f^n)\cap K}\exp
      S_n\varphi^u(x)\right).
  \end{equation}
  Since $K$ is hyperbolic, the potential
  $\varphi^u_K=-\log|\det Df|E^u|\colon K\to\bR$ is H\"older continuous.
  Therefore, Proposition~\ref{ha} and~\eqref{fix2} imply that
  \begin{equation*}
    P_{\rm SP} (\varphi^u,\alpha,c)\leq P_{\rm top}(f|K,\varphi^u_K).
  \end{equation*}
  On the other hand, the variational principle gives
  \begin{equation*}
    P_{\rm top}(f|K,\varphi^u_K)=
    \sup_{\nu}\left(h_\nu(f)+\int\varphi^ud\nu\right),
  \end{equation*}
  where the supremum is taken over all $\nu\in\cM_{\rm E}$ which are supported on
  $K$ (and which, in particular, are of saddle type).  Moreover, by construction of $K$ we have that $\alpha\leq\chi(\nu)$ holds for these measures.
  We conclude that
  \begin{equation*}
  \lim_{c\to 0} P_{\rm SP} (\varphi^u,\alpha,c) \le \sup_{\nu}\left(h_\nu(f)+\int\varphi^ud\nu\right),
  \end{equation*}
  where the supremum is taken over all
  measures $\nu\in\cM_{\rm E}$ of saddle type with $\alpha\le\chi(\nu)$.
\end{proof}

\ni
Under the assumptions of Theorem~\ref{varprinc} we call
\begin{equation*}
P_{\rm SP} (\varphi^u)\eqdef\lim_{\alpha\to 0}\lim_{c\to 0} P_{\rm SP} (\varphi^u,\alpha,c)
\end{equation*}
the \emph{volume pressure} of $f$.\\[0.2cm]
{\it Remark. }
We note that if $\Lambda$ is a locally maximal hyperbolic set such that
$f|\Lambda$ is topologically mixing  then the volume pressure coincides with
the topological pressure of the potential $\varphi^u_\Lambda=-\log|\det
Df|E^u|\colon \Lambda\to\bR$.\\[0.2cm]
As an immediate consequence of the proof of Theorem~\ref{varprinc} we obtain
the following ``inverse" variational inequality for the volume pressure.

\begin{corollary}\label{cor2}
  Let $f$ and $\Lambda$ be as in Theorem~\ref{varprinc}, and assume that $f$
  has a saddle point in $\Lambda$. Then
  \begin{equation*}
    P_{\rm SP} (\varphi^u)
    \le\sup_{K} P_{\rm top}(f|K,-\log\lvert\det Df|E^u\rvert),
  \end{equation*}
  where the supremum is taken over all compact
  hyperbolic sets $K\subset \Lambda$.
\end{corollary}

Corollary~\ref{cor2} immediately implies the following.

\begin{corollary}\label{cor3}
  Let $f$ and $\Lambda$ be as in Theorem~\ref{varprinc}, and assume that $f$
  has a saddle point in $\Lambda$. Then
  \begin{equation*}
    P_{\rm SP} (\varphi^u)
    \leq\sup_\nu \left( h_\nu(f)+\int_\Lambda \varphi^ud\nu\right),
  \end{equation*}
  where the supremum is taken over all
  measures $\nu\in\cM_{\rm E}$ of saddle type.
\end{corollary}

\section{Rate of escape from neighborhoods and dimension of invariant
  sets}\label{sec:5}

In this section we discuss relations between the attraction properties of the
invariant set $\Lambda$ and the dimension of certain subsets of $\Lambda$.

Let $\Lambda\subset M$ be a compact locally maximal  $f$-invariant set and let
$U\subset M$ be an open neighborhood of $\Lambda$ such that $\Lambda =
\bigcap_{n\in\bZ}f^n(U)$.
Given an open neighborhood $V\subset U$ of $\Lambda$ we define the \emph{upper
  (exponential) escape rate} from $V$ by
\begin{equation}\label{defvol1}
\overline{E}(V) \eqdef
\limsup_{n\to\infty}\frac{1}{n}\log
\left(\vol\bigcap_{k=0}^{n-1}f^{-k}(V)\right),
\end{equation}
where $\vol$ denotes the volume induced by the Riemannian metric on
$M$. Analogously, we define the lower escape rate $\underline{E}(V)$ by
replacing the limes superior in~\eqref{defvol1} with the limes inferior.  By
definition, we have that $\underline{E}(V)\leq\overline{E}(V)\leq
0$. Under the assumption that the upper escape rate is
strictly negative we obtain a non-trivial upper bound for the  upper box dimension
of $\Lambda$ (see Theorem~\ref{thSW} below). Define
\begin{equation}\label{defs}
  s\eqdef
  \lim_{n\to\infty} \frac{1}{n}\log\left(\max_{x\in \Lambda} \lVert
    Df^n(x)\rVert\right).
\end{equation}
Note that $s$ is well-defined. This follows from the sub-additivity of the
sequence $(\varphi_n)_n$ given by $\varphi_n=\log\left(\max_{x\in\Lambda}\lVert
Df^n(x)\rVert\right)$ (see e.g.~\cite{Wal:81}).

\begin{theorem}\label{thSW}
  Let $f\colon M\to M$ be a $C^2$-diffeomorphism, and let
  $\Lambda\subset M$ be a  compact locally maximal  $f$-invariant set
  containing a point with a positive Lyapunov exponent. Then
  \begin{equation}\label{maininequality}
    \overline{\dim}_{\rm B} \Lambda \leq \dim M + \frac{\overline{E}(V)}{s}.
  \end{equation}
  In particular, if $\overline{E}(V)<0$ then $\overline{\dim}_{\rm B}
  \Lambda<\dim M$.
\end{theorem}

\begin{proof}
  First,  we note that since $\Lambda$ contains a point with a positive
  Lyapunov exponent, it follows that $s>0$. Let $\delta>0$. By a simple
  continuity argument there exist $\epsilon>0$ and $n(\delta)\in \bN$ such
  that for all $y\in B(\Lambda,\epsilon) \eqdef \bigcup_{x\in \Lambda}
  B(x,\epsilon) $ we have
  \begin{equation*}
    \lVert Df^{n(\delta)}(y)\rVert < \exp(n(\delta)(s+\delta)).
  \end{equation*}
  From now on we consider the map $g=f^{n(\delta)}$. Note that $\Lambda$
  is also a compact invariant set of $g$.
  It follows from~\eqref{defvol1} and~\eqref{defs} that $\overline{E}_g(V)\leq
  n(\delta)\overline{E}_f(V)$ and $s_g=n(\delta) s_f$. Therefore, it suffices
  to prove~\eqref{maininequality} for $g$. We continue to use the notation $s$
  and $\overline{E}(V)$ for $g$ instead of $f$.

  Set $V=B(\Lambda,\epsilon)$. By making $\epsilon$ smaller if necessary we can
  assure that $V\subset U$. It follows from the definition of $\overline{E}(V)$
  that  if $n$ is sufficiently large then
  \begin{equation}\label{eqbo4}
    \vol\left(\bigcap_{k=0}^{n-1} f^{-k}(V)\right)<
    \exp(n(\overline{E}(V)+\delta)).
  \end{equation}
  For $n\in \bN$ we define real numbers
  \[
  r_n = \frac{\epsilon}{\exp(n(s+\delta))}
  \]
  and neighborhoods $B_n = B(\Lambda,r_n)$ of $\Lambda$. Let $y\in B_n$. Then
  there exists $ x\in\Lambda$ with $d(x,y) < r_n$. An elementary induction
  argument in combination with the mean-value theorem implies
  $d(g^{i}(x),g^{i}(y)) < \epsilon$ for all $i\in\{0,\dots,n-1\}$. Hence
  $B_n\subset \bigcap_{k=0}^{n-1} g^{-k}(V)$, and~\eqref{eqbo4} implies that
  \begin{equation*}
    \vol(B_n)<
    \exp\left(n\left(\overline{E}(V)+\delta\right)\right)
  \end{equation*}
  for sufficiently large $n$.
  Let us recall that for $t\in [0,\dim M]$ the  $t$-dimensional upper Minkowski
  content of a relatively compact set $A\subset M$ is defined by
  \[
  M^{*t}(A)= \limsup_{\rho\to 0} \frac{\vol(A_\rho)}{(2\rho)^{\dim M-t}},
  \]
  where $A_\rho=\{y\in M\colon\exists x\in A\colon d(x,y)\leq \rho\}$. Let
  $t\in [0,\dim M]$ and $\rho_n={\textstyle\frac{r_n}{2}}$ for all
  $n\in\bN$. Then we have
  \begin{equation}\label{glej-4}
    \begin{split}
      &M^{*t} (\Lambda)
      =\limsup_{\rho\to 0}\frac{\vol(\Lambda_\rho)}{(2\rho)^{\dim M-t}}\\
      &\leq \limsup_{n\to\infty }
      \frac{\vol(\Lambda_{\rho_n})}{(2\rho_{n+1})^{\dim M-t}}
      \leq
      \limsup_{n\to \infty} \frac{\vol(B_n)}{(r_{n+1})^{\dim M-t}}\\
      &\leq  \frac{\exp((\dim M-t)(s+\delta))}{\epsilon^{\dim M-t}}
      \lim_{n\to\infty}\left(\exp\left((\dim M-t)(s+\delta) +
          \overline{E}(V)+\delta \right)\right)^n.
    \end{split}
  \end{equation}
  Assume  $t>\dim M + (\overline{E}(V)+\delta)/(s+\delta)$.  Then
  \[
  \exp\left((\dim M-t)(s+\delta) + \overline{E}(V)+\delta \right) < 1,
  \]
  which implies $M^{*t}(\Lambda) = 0$. Hence $t\geq \overline{\dim}_{\rm B}
  \Lambda$ (see~\cite{Ma}). Finally, the fact that $\delta$ was arbitrary
  completes the proof.
\end{proof}
\noindent
{\it Remarks. }\\
(i) In the proof of Theorem~\ref{thSW} we have used similar ideas as
in~\cite{ShaWol:05} in the context of hyperbolic sets.\\
(ii) It is easy to see that the analog of Theorem~\ref{thSW} holds for the
lower box dimension  with $\overline{E}(V)$ replaced by $\underline{E}(V)$ in
\eqref{maininequality}.\\
(iii) The dimension estimate in Theorem~\ref{thSW} remains true when $f$ is locally a $C^1$-diffeomorphism.
\\[0.2cm]

It is a result of Young~\cite[Theorem 4 (1)]{You:90}  that
\begin{equation}\label{you}
\sup_{\nu\in\cM_{\rm E}}\left(h_\nu(f)+\int_\Lambda\varphi^ud\nu\right)\le
\underline E(V) \le \overline E(V) \le 0.
\end{equation}
Therefore, Theorem~\ref{varprinc} implies that if $P_{\rm SP}(\varphi^u)=0$ then
$\underline E(V)= \overline E(V)=0$ and thus, the
volume of $\bigcap_{k=0}^{n-1} f^{-k}(V)$ is shrinking at a  sub-exponential
rate. 
It is a classical result of Bowen~\cite{B} that in the case when
$\Lambda$ is a hyperbolic set of $f$ then
\begin{equation}\label{youu}
\sup_{\nu\in\cM_{\rm E}}\left(h_\nu(f)+\int_\Lambda\varphi^ud\nu\right)=
\underline E(V) = \overline E(V),
\end{equation}
and  thus the volume of $\bigcap_{k=0}^{n-1} f^{-k}(V)$ shrinks at an
exponential rate (in which case $\Lambda$ is an attractor) if and only if
$P_{\rm SP}(\varphi^u)=0$. 

Baladi et al.~\cite{BalBonSch:99} give an example of a compact locally maximal
$f$-invariant set $\Lambda$ of a $C^\infty$-surface diffeomorphism and an
invariant measure $\mu$, supported on $\Lambda$, which attains the supremum
in~\eqref{you} with the property that $h_\mu(f)+\int_\Lambda\varphi^ud\mu
<\overline E(V)=0$ for an arbitrarily small neighborhood $V$ of $\Lambda$
(another example is given by the ``figure-8 attractor",
see~\cite[p.~140]{Kat:80}). Obviously, in such
a situation Theorem~\ref{thSW} does not provide a nontrivial upper bound for
the box dimension of $\Lambda$. 
The following example shows that in general such a nontrivial upper bound does
not exist.

\begin{example}\label{ex:3}{\rm 
    Let $M$ be a smooth compact surface and let $f:M\to M$
    be a $C^2$ diffeomorphism whose non-wandering set is a horseshoe
    $K$. Define $\Lambda=M$. Then  $\overline{\dim}_{\rm B} \Lambda= \dim M =
    2$ and $\overline E(V)=0$, however $P_{\rm SP}(\varphi^u)= P_{\rm
      top}(f|K,\varphi^u)<0$.
}\end{example} 

We end the paper with the following \\[0.2cm]
\noindent {\bf Open problem. }Does  $P_{\rm SP}(\varphi^u)<0$
imply $\overline{\dim}_{\rm B}
  \overline{\SPer(f)}<\dim M$?\\[0.3cm]
                                
\ni
{\bf Acknowledgment. } We would like to thank the referee for several
comments and improvements.


\begin{thebibliography}{99}
\bibitem{BalBonSch:99} V.~Baladi, Ch.~Bonatti, and B.~Schmitt, \emph{Abnormal
    escape rates from nonuniformly hyperbolic sets}, Ergodic Theory
  Dynam. Systems \textbf{19} (1999), 1111--1125.
\bibitem{Bar:96} L.~Barreira, \emph{A non-additive thermodynamic formalism and
    applications to dimension theory of hyperbolic dynamical systems},
  Ergodic Theory Dynam. Systems \textbf{16} (1996), 871--927.
\bibitem{BarPes:05} L.~Barreira and Y.~Pesin, \emph{Smooth ergodic theory and
    nonuniformly hyperbolic dynamics}, in Handbook of Dynamical Systems
  \textbf{1B}, B.~Hasselblatt and A.~Katok eds.,  Elsevier, 2006.
\bibitem{BonDiaPuj:03} C.~Bonatti, L.~D\'iaz, and E.~Pujals, \emph{A
    $C^1$-generic dichotomy for diffeomorphisms: weak forms of hyperbolicity
    or infinitely many sinks of sources}, Ann. Math. \textbf{158} (2003),
  355--418.
\bibitem{B} R.~Bowen, \emph{Equilibrium states and the ergodic theory
    of Anosov diffeomorphisms}, Lecture Notes in Mathematics \textbf{470},
  Springer, 1975.
\bibitem{BR} R. Bowen and D. Ruelle, \emph{The ergodic theory of Axiom A
    flows}, Invent. Math. \textbf{29} (1975), 181--202.
\bibitem{Chu} Y. M. Chung, \emph{Expanding periodic orbits with small
    exponents}, J. Difference Equ. Appl. \textbf{9} (2003), 337--341.
\bibitem{ChuHir:03}  Y.~M.~Chung and M.~Hirayama,
  \emph{Topological entropy and periodic orbits of saddle type for surface
    diffeomorphisms}, Hiroshima Math. J. \textbf{33} (2003), 189--195.
\bibitem{Kat:80} A.~Katok, \emph{Lyapunov exponents, entropy and periodic
    orbits for diffeomorphisms}, Publ. Math., Inst. Hautes
  \'Etud. Sci. \textbf{51} (1980), 137--173.
\bibitem{KatHas:95} A.~Katok and B.~Hasselblatt, \emph{Introduction to the
    Modern Theory of Dynamical Systems}, Encyclopedia of Mathematics
  and Its Applications \textbf{54}, Cambridge University Press, 1995.
\bibitem{Ma} P.~Mattila, \emph{Geometry of sets and measures in
    Euclidean spaces. Fractals and rectifiability}, Cambridge
  University Press, 1995.
\bibitem{P} Y.~Pesin, \emph{Dimension theory in dynamical systems:
    Contemporary Views and applications}, Lectures in Mathematics, Chicago
  University Press, 1997. 
\bibitem{ShaWol:05} R.~Shafikov and Ch.~Wolf, \emph{Stable sets, hyperbolicity
    and dimension}, Discrete Contin. Dynam. Systems \textbf{12} (2005),
  403--412.
\bibitem{Wal:81} P.~Walters, \emph{An introduction to ergodic theory},
  Graduate Texts in Mathematics \textbf{79}, Springer, 1981.
\bibitem{You:90} L.-S.~Young, \emph{Large deviations in dynamical systems},
  Trans. Amer. Math. Soc. \textbf{318} (1990), 525--543.
\end{thebibliography}
\end{document}